\documentclass[12pt]{article}

\usepackage{amssymb}

\makeatletter
 
 \@addtoreset{equation}{section}
\makeatother

\newcommand{\mod}{\mbox{\rm mod}\,}

\newcommand{\li}{\mbox{\rm li}\,}

\newcommand{\Gal}{\mbox{\rm Gal}}
\newcommand{\Aut}{\mbox{\rm Aut}}
\newcommand{\id}{\mbox{\rm id}}
\newcommand{\Frac}{\displaystyle\frac}
\newcommand{\psum}{\mathop{{\sum}'}}
\renewcommand{\pmod}[1]{%
\ (\mbox{\rm mod}\ #1)}

\newcommand{\prf}{\noindent{\bf Proof. }}

\newcommand{\rem}{\noindent{\bf Remark. }}

\newcommand{\qed}{\hbox{\rule[-2pt]{5pt}{11pt}}}

\newtheorem{dfn}{Definition}[section]
\newtheorem{thm}[dfn]{Theorem}
\newtheorem{prop}[dfn]{Proposition}
\newtheorem{cor}[dfn]{Corollary}
\newtheorem{lem}[dfn]{Lemma}

\newtheorem{hyp}[dfn]{Hypothesis}

\begin{document}
\title{On a distribution property of the residual order of $a\pmod p$}
\author{Koji Chinen\footnotemark[1] and Leo Murata\footnotemark[7]}
\date{}
\maketitle
\begin{abstract}
Let $a$ be a positive integer with $a\ne1$ and $Q_a(x;k,l)$ be the set of primes $p\leq x$ such that the residual order of $a$ in ${\bf Z}/p{\bf Z}^\times$ is congruent to $l$ $\bmod$ $k$. It seems that no one has ever considered the density of $Q_a(x;k,l)$ for $l\ne0$ when $k\geq 3$. In this paper, the natural densities of $Q_a(x;4,l)$ ($l=0,1,2,3$) are considered. 
We assume $a$ is square free and $a\equiv 1\pmod 4$. Then, for $l=0,2$ , we can prove unconditionally that their natural densities are equal to 1/3. On the contrary, for $l=1,3$ , we assume Generalized Riemann Hypothesis, then we can prove their densities are equal to 1/6.
\end{abstract}
\footnotetext[1]{Department of Mathematics, Faculty of Engineering, Osaka Institute of Technology. Omiya, Asahi-ku, Osaka 535-8585, Japan. E-mail: YHK03302@nifty.ne.jp}
\footnotetext[7]{Department of Mathematics, Faculty of Economics, Meiji Gakuin University, 1-2-37 Shirokanedai, Minato-ku, Tokyo 108-8636, Japan. E-mail: leo@eco.meijigakuin.ac.jp}

This manuscript is the one which we submitted to Crelle Journal in March 2001. The second author talked on this subject at Oberwolfach "Theory of the Riemann Zeta and Allied Functions" at 20.09.2001.
\section{Introduction}
Let ${\bf P}$ be the set of all prime numbers. 

For a fixed natural number $a\geq2$, we can define two functions, $I_a$ and $D_a$, from ${\bf P}$ to ${\bf N}$:
\begin{equation}\begin{array}{rrrl}
I_a: & p & \mapsto & I_a(p)=\left|({\bf Z}/p{\bf Z})^\times:\langle a \pmod p\rangle\right|\\[2mm]
     &   &         & \mbox{(the residual index mod $p$ of $a$)},\\[3mm]
D_a: & p & \mapsto & D_a(p)=\sharp\langle a\pmod p\rangle\\[2mm]
     &   &         & \mbox{(the order of the class $a$ ($\mod p$) in $({\bf Z}/p{\bf Z})^\times$)},
\end{array}
\end{equation}
where $({\bf Z}/p{\bf Z})^\times$ denotes the set of all invertible residue classes  $\mod p$, and $|\ :\ |$ the index of the subset. 

We have a simple relation 
\begin{equation}\label{eq:DaIa=p-1}
I_a(p)D_a(p)=p-1,
\end{equation}
but both of these functions fluctuate quite irregularly. C. F. Gauss already noticed that $I_{10}(p)=1$ happens rather frequently. And the famous Artin's conjecture for primitive roots asks whether the cardinality of the set
\begin{equation}\label{eq:Na(x)}
N_a(x):=\bigl\{p\leq x\ ;\ I_a(p)=1 \bigr\}
\end{equation}
tends to $\infty$ or not as $x\to\infty$. On the assumption of the Generalized Riemann Hypothesis for a certain type of Dedekind zeta functions, C. Hooley \cite{Ho} succeeded in calculating the natural density of $N_a(x)$. There are various variations of Artin's Conjecture, among which two papers Lenstra \cite{Le} and Murata \cite{Mu} considered the surjectivity of the map $I_a$. For any natural number $n$, we define 
\begin{equation}\label{eq:Na(x;n)}
N_a(x;n):=\bigl\{p\leq x\ ;\ I_a(p)=n \bigr\}.
\end{equation}
Then their results show that, for a square free $a$ with $a\not\equiv 1\pmod 4$, we have, under GRH, an asymptotic formula
\begin{equation}\label{eq:Le-Mu}
\sharp N_a(x;n)\sim C_a^{(n)}\li x
\end{equation}
and $C_a^{(n)}>0$, where $\li x:=\int_2^x(\log t)^{-1}dt$ and the constant $C_a^{(n)}$ depends on $a$ and $n$. Therefore, for such an $a$, the map $I_a$ is surjective from ${\bf P}$ onto ${\bf N}$. 

And the surjectivity of the map $D_a$ is also proved by many authors. They proved that, except for at most finitely many $n$'s, the map $D_a$ is surjective from ${\bf P}$ onto ${\bf N}$. 

Thus these two maps are surjective for those $a$'s, but between their surjective-properties we notice a big difference. Under GRH, for any $n\in{\bf N}$, (\ref{eq:Le-Mu}) means that 
\begin{equation}\label{eq:I_a_inv}
I_a^{-1}(n)=\bigl\{p\in{\bf P}\ ;\ I_a(p)=n\bigr\}
\end{equation}
contains infinite elements, but on the contrary, the set 
\begin{equation}\label{eq:D_a_inv}
D_a^{-1}(n)=\bigl\{p\in{\bf P}\ ;\ D_a(p)=n\bigr\}
\end{equation}
contains only a finite number of elements. In fact, if $D_a(p)=n$, then
$$n+1\leq p\leq a^n.$$
And recent study on cryptography shows that characterizing $D_a$ is very difficult. 

For the purpose of considering the distribution property of the map $D_a$, here we take an arbitrary natural number $k\geq 2$ and an arbitrary residual class $l\pmod k$ and consider the asymptotic behavior of the cardinality of the following set: 
\begin{equation}\label{eq:Q_a_def} 
Q_a(x;k,l):=\bigl\{p\leq x\ ;\ D_a(p)\equiv l\pmod k \bigr\}.
\end{equation}
It is more than 50 years ago, W. Sierpinski first considered this problem and H. Hasse proved, by our notations, that, for an odd prime $q$, 
$$\mbox{the Dirichlet density of }Q_a(x;q,0)=\frac{q}{q^2-1}$$
(\cite{Ha1} and \cite{Ha2}). Odoni \cite{Od} proved the existence of the natural density of $Q_a(x;q,0)$, and he obtained a similar results on $Q_a(x;k,0)$ for a composite square free moduli $k$. 

In this paper we take $k=4$ and consider the distribution property of $Q_a(x;4,l)$ for all residue classes $l=0,1,2,3$. Our results consist of two theorems:
\begin{thm}\label{th:uncond}
We assume $a$ is a square free positive integer with $a\geq 3$. When $l=0,2$, we have 
$$\sharp Q_a(x;4,l)=\frac{1}{3}\li x+O\left(\frac{x}{\log x\log\log x}\right).$$
\end{thm}
\begin{thm}\label{th:cond}
Let $a$ be as above. We assume GRH and further assume $a\equiv 1\pmod 4$. Then, for $l=1,3$, we have 
$$\sharp Q_a(x;4,l)=\frac{1}{6}\li x+O\left(\frac{x}{\log x\log\log x}\right).$$
\end{thm}
Here GRH means:
\begin{hyp}[Generalized Riemann Hypothesis]\label{th:GRH}
For any positive integers $k$ and $m$ with $k|m$, we assume that the Riemann Hypothesis holds for the Dedekind zeta function $\zeta_K(s)$ for the field $K={\bf Q}(\zeta_m, a^{1/k})$ where $\zeta_m=\exp(2\pi i/m)$. 
\end{hyp}
As we mentioned above, Hasse and Odoni investigated $\sharp Q_a(x;q,l)$ with $l=0$, and for $l\ne 0$, the distribution property of $\sharp Q_a(x;q,l)$ remains unknown so far. 

When $l=0$, the condition ``$D_a(p)\equiv 0\pmod q$'' can be reformulated in the notation of algebraic number theory without much difficulty. In fact, for a prime $p$ with $q||p-1$, using the relation (\ref{eq:DaIa=p-1}),
$$D_a(p)\equiv 0\pmod q \quad \Leftrightarrow \quad q\nmid I_a(p),$$
and we can count the number of such primes $p$'s by
$$\sharp\bigl\{p\leq x\ ;\ q||p-1 \bigr\}-
\sharp\bigl\{p\leq x\ ;\ q||p-1,\ q|I_a(p) \bigr\}.$$
The last condition $q|I_a(p)$ means that $a$ is a $q$-th power residue modulo $p$, so we can utilize the prime ideal theorem etc. 

On the contrary, when $l\ne0$, the reformulation of ``$D_a(p)\equiv l\pmod q$'' needs a rather complicated procedure (see Lemma \ref{th:Q_a_sieve} and compare (i) and (iii) ). Moreover, e.g. the calculation of each $\sharp N_a(x; 2^f+l\cdot 2^{f+2}; 1+2^f\pmod{2^{f+2}})$ which appears in (\ref{eq:sieve_4_1}) requires a consideration on the generalized Artin's conjecture $N_a(x;2^f+l\cdot 2^{f+2})$ in the special residue class $p\equiv 1+2^f\pmod{2^{f+2}}$. That is why we need GRH in Theorem \ref{th:cond}. It seems an interesting phenomenon that, after all, we arrived at such a simple result as Theorem \ref{th:cond}. 

And our result as well as Odoni's result shows that the value distribution of the map $D_a$ is rather irregular. 

We prepare some preliminary lemmas in Section 2, prove Theorem \ref{th:uncond} in Section 3, and prove the conditional result Theorem \ref{th:cond} in Section 4. In Section 5, we mention some numerical examples. 
For our results, see also [ 9 ].
\section{Preliminaries}
Throughout this paper we fix a square free integer $a\geq 3$, and  $p$ denotes an odd prime which does not divide $a$. For $k\in{\bf Z}$, let $\zeta_k=\exp(2\pi i/k)$. We denote Euler's totient and the M\"obius function by $\varphi(k)$ and $\mu(k)$, respectively. For a prime power $q^e$, $q^e||m$ means that $q^e|m$ and $q^{e+1}\nmid m$. 

Let $K$ be an algebraic number field. Then we define 
\begin{equation}\label{eq:pi_x_K}
\pi(x,K)=\sharp\bigl\{\mathfrak p: \mbox{ a prime ideal in }K,\ N\mathfrak p\leq x\bigr\}
\end{equation}
and 
\begin{equation}\label{eq:pi1_x_K}
\pi^{(1)}(x,K)=\sharp
\bigl\{\mathfrak p: \mbox{ a prime ideal of degree 1 in }K,\ N\mathfrak p\leq x\bigr\}
\end{equation}
where $N\mathfrak p$ is the (absolute) norm of $\mathfrak p$. Moreover let $L/K$ be a finite Galois extension. Then for a prime ideal $\mathfrak p$ in $K$, we define the Frobenius symbol by
\begin{equation}\label{eq:frob}
({\mathfrak p}, L/K)=\left\{
     \begin{array}{ll}\sigma\in\Gal(L/K)\ ; & {\mathfrak q}^\sigma={\mathfrak q} \mbox{ for some prime } {\mathfrak q} \mbox{ in }L \mbox{ above }{\mathfrak p},\\
\ & \alpha^\sigma \equiv \alpha^{N{\mathfrak p}} \pmod{\mathfrak q} \mbox{ for all }\alpha\in L
     \end{array}\right\}.
\end{equation}
This notation is due to Lenstra \cite{Le}. 

Next we introduce some preliminary results. In the course of our proof, we need the exact value of the extension degree of a certain type of Kummer fields. In the following lemma only, we include the case $a$ is not a square free integer. 
\begin{lem}\label{th:moree_lemma}
Let $k,r\in{\bf N}$ with $k|r$. We assume $a$ is not a perfect $h$-th power with $h\geq 2$. And $a_1$ being the square free part of $a$ (i.e. $a=a_1a_2^2$ with $a_1$: square free), we put
$$h_1=\left\{\begin{array}{ll}
  2a_1, & \mbox{if }a_1\equiv 1\pmod 4,\\
  4a_1, & \mbox{otherwise.}
\end{array}\right.$$
Then we have
\begin{equation}\label{eq:moree_lemma}
[{\bf Q}(\zeta_r,a^{1/k}):{\bf Q}]=\left\{\begin{array}{ll}
k\varphi(r), & \mbox{if $k$ is odd,}\\[2mm]
k\varphi(r), & \mbox{if $k$ is even and }h_1\nmid r,\\[2mm]
k\varphi(r)/2, & \mbox{if $k$ is even and }h_1|r.
\end{array}\right.
\end{equation}
\end{lem}
\prf See Moree \cite[Lemma 2]{Mo} or Murata \cite[Section 3]{Mu}. \qed

\medskip
\begin{thm}\label{th:PIT}
For a prime $q$ and $i,j\in{\bf N}\cup \{0\}$, we define an extension field
$$K_{i,j}^{(q)}={\bf Q}(\zeta_{q^i}, \zeta_{q^j}, a^{1/q^j}),$$
and we put
\begin{eqnarray*}
n&=&[K_{i,j}^{(q)}:{\bf Q}],\\
D&=&\mbox{the discriminant of }K_{i,j}^{(q)}.
\end{eqnarray*}
Then, under the condition 
$$x\geq\exp\bigl(10n\log^2|D|\bigr),$$
we have
$$\pi^{(1)}(x,K_{i,j}^{(q)})=\li x+O\bigl(nxe^{-c\sqrt{\log x}/n^2}\bigr),$$
where the constant implied by $O$-symbol and the positive constant $c$ depend only on $a$ and $q$.
\end{thm}
\prf For the field $K_{i,j}^{(q)}$, we have an estimate
$$|D|\leq(n^2|a|)^n.$$ 
Then Theorems 1.3 and 1.4 of Lagarias-Odlyzko \cite{LaOd} give the desired formula. \qed

\medskip
And we need the Chebotarev density theorem with GRH: 
\begin{thm}[Chebotarev density theorem, GRH]\label{th:Chebotarev}
Let $K$ be an algebraic number field, $L/K$ be a finite Galois extension and $C$ be a conjugacy class in $G=\Gal(L/K)$. We define $\pi(x;L/K,C)$ by
\begin{equation}
\pi(x;L/K,C)=\sharp\{{\mathfrak p}: \mbox{ a prime ideal in $K$, unramified in $L$, }({\mathfrak p}, L/K)=C,\ N{\mathfrak p}\leq x\}.
\label{eq:def_pi_C}
\end{equation}
Then, under GRH for the field $L$, we have
\begin{equation}
\pi(x;L/K,C)=\frac{\sharp C}{\sharp G}\li x+O\left(\frac{\sharp C}{\sharp G}\sqrt x\log(d_Lx^{n_L})+\log d_L\right),\quad\mbox{as }x\to\infty,
\label{eq:Chebotarev}
\end{equation}
where $d_L$ is the discriminant of $L$ and $n_L=[L:{\bf Q}]$. 
\end{thm}
\prf Lagarias-Odlyzko \cite[Theorem 1.1]{LaOd}. \qed

\medskip
Here we recall the set of primes $N_a(x;n)$ which we defined in (\ref{eq:Na(x;n)}). 
\begin{lem}\label{th:lemma_I_a(p)}
We assume GRH. Let $a$ be a square free integer $\geq2$, $\psi(x)$ be a monotone increasing positive function which satisfies
$$\lim_{x\to\infty}\psi(x)=+\infty\quad\mbox{ and }\quad
\psi(x)\ll(\log x)^{\frac14}.$$
Then we have 
$$\sharp\bigl\{p\leq x\ ;\ I_a(p)\geq \psi(x)\bigr\}\ll\frac{\pi(x)}{\psi(x)},$$
where the constant implied by $\ll$-symbol is absolute. 
\end{lem}
\prf Let $y$ be the largest integer not exceeding $\psi(x)$. We have 
\begin{equation}\label{eq:I_a(p)_and_y}
\bigl\{p\leq x\ ;\ I_a(p)\geq y\bigr\}
=\bigl\{p\leq x\ ;\ p\nmid a\bigr\}-\bigcup_{n=1}^{y-1}N_a(x;n),
\end{equation}
where $\cup_{n=1}^{y-1}$ is a disjoint union. Then Theorems 1 and 2 of Murata \cite{Mu} prove, with $\varepsilon=1$, 
\begin{equation}\label{eq:Murata_Th1}
\sharp N_a(x;n)
=C_a^{(n)}\li x+O\left(\{n\log\log x+\log a\}\frac{x}{\log^2 x}\right)
\end{equation}
and 
\begin{equation}\label{eq:Murata_Th2}
\sum_{n\leq y}C_a^{(n)}=1+O\left(\frac{1}{y}\right).
\end{equation}
Thus from (\ref{eq:I_a(p)_and_y}), we have
\begin{eqnarray*}
\sharp\bigl\{p\leq x\ ;\ I_a(p)\geq \psi(x)\bigr\}
&=&\pi(x)-\left(1+O\left(\frac{1}{y}\right)\right)\li x
  +O\left(\frac{x\log\log x}{\log^2 x}\cdot\sum_{n\leq y}n\right)\\
&=&O\left(\frac{1}{y}\pi(x)\right)
  +O\left(\frac{x\log\log x}{(\log x)^{3/2}}\right)
  =O\left(\frac{\pi(x)}{\psi(x)}\right).\ \qed
\end{eqnarray*}
\section{Proof of Theorem \ref{th:uncond}}
Generally speaking, the condition ``$D_a(p)\equiv j\pmod 4$'' is rather difficult to handle. So, using the relation (\ref{eq:DaIa=p-1}), we transform the condition on $D_a(p)$ into some conditions on $I_a(p)$. Here we introduce the set
\begin{equation}\label{eq:Na(xnst)}
N_a(x;n;s\pmod t):=\bigl\{p\leq x\ ;\ p\in N_a(x;n),\ p\equiv s\pmod t \bigr\},
\end{equation}
which is a generalization of the set $N_a(x;n)$ which appeared in Artin's conjecture for primitive roots. 

We can prove the following lemma, which is the starting point of our proof:
\begin{lem}\label{th:Q_a_sieve}
For any $x>0$ we have\\
(i)\begin{eqnarray}
\sharp Q_a(x;4,0)&=&\sharp\bigl\{p\leq x\ ;\ p\equiv 1\pmod 4\bigr\}\nonumber\\
 &-& \sum_{j\geq1}
   \sharp\bigl\{p\leq x\ ;\ p\equiv 1\pmod{2^{j+1}},
   \ 2^j|I_a(p)\bigr\}\nonumber\\
 &+&\sum_{j\geq1}
   \sharp\bigl\{p\leq x\ ;\ p\equiv 1\pmod{2^{j+2}},\ 2^j|I_a(p)\bigr\},
\label{eq:sieve_4_0}
\end{eqnarray}
(ii)\begin{eqnarray}
\sharp Q_a(x;4,2)
 &=&\sum_{j\geq1}
   \sharp\bigl\{p\leq x\ ;\ p\equiv 1\pmod{2^j},
   \ 2^{j-1}|I_a(p)\bigr\}\nonumber\\
 &-& \sum_{j\geq1}
   \sharp\bigl\{p\leq x\ ;\ p\equiv 1\pmod{2^{j+1}},
   \ 2^{j-1}|I_a(p)\bigr\}\nonumber\\
 &-& \sum_{j\geq1}
   \sharp\bigl\{p\leq x\ ;\ p\equiv 1\pmod{2^j},\ 2^j|I_a(p)\bigr\}\nonumber\\
 &+&\sum_{j\geq1}
   \sharp\bigl\{p\leq x\ ;\ p\equiv 1\pmod{2^{j+1}},\ 2^j|I_a(p)\bigr\},
\label{eq:sieve_4_2}
\end{eqnarray}
(iii)\begin{eqnarray}
\sharp Q_a(x;4,1)&=&\sum_{f\geq1}\sum_{l\geq0}\sharp 
  N_a(x;2^f+l\cdot 2^{f+2};1+2^f\pmod{2^{f+2}})\nonumber\\
&+&\sum_{f\geq1}\sum_{l\geq0}\sharp 
  N_a(x;3\cdot 2^f+l\cdot 2^{f+2};1+3\cdot 2^f\pmod{2^{f+2}}),
\label{eq:sieve_4_1}
\end{eqnarray}
(iv)\begin{eqnarray}
\sharp Q_a(x;4,3)&=&\sum_{f\geq1}\sum_{l\geq0}\sharp 
  N_a(x;3\cdot 2^f+l\cdot 2^{f+2};1+2^f\pmod{2^{f+2}})\nonumber\\
&+&\sum_{f\geq1}\sum_{l\geq0}\sharp 
  N_a(x;2^f+l\cdot 2^{f+2};1+3\cdot 2^f\pmod{2^{f+2}}).
\label{eq:sieve_4_3}
\end{eqnarray}
\end{lem}
\prf We can prove these formulas in a similar manner. Here we show only the proof of (i). From the condition 
$$D_a(p)\equiv 0\pmod 4,$$
it is necessary that $p\equiv 1\pmod 4$. So we consider a prime such that $2^j||p-1$, $j\geq 2$. Then, with the relation (\ref{eq:DaIa=p-1}), we have
$$D_a(p)\equiv 0\pmod 4 \quad \Leftrightarrow \quad 2^{j-1}\nmid I_a(p).$$
Hence we have
\begin{eqnarray*}
Q_a(x;4,0)&=&\bigcup_{j\geq 2}
\bigl\{p\leq x\ ;\ 2^j||p-1,\ 2^{j-1}\nmid I_a(p)\bigr\}\\
  &=&\bigcup_{j\geq 2}\Bigl(\bigl\{p\leq x\ ;\ p\equiv 1\pmod{2^j},\ 
    2^{j-1}\nmid I_a(p)\bigr\}\\
  & &-\bigl\{p\leq x\ ;\ p\equiv 1\pmod{2^{j+1}},\ 
    2^{j-1}\nmid I_a(p)\bigr\}\Bigr)\\
 & & \\
&=&\bigl\{p\leq x\ ;\ p\equiv 1\pmod 4\bigr\}\\
 &-& \bigcup_{j\geq 2}
   \bigl\{p\leq x\ ;\ p\equiv 1\pmod{2^{j}},\ 2^{j-1}|I_a(p)\bigr\}\\
 &+&\bigcup_{j\geq 2}
   \bigl\{p\leq x\ ;\ p\equiv 1\pmod{2^{j+1}},\ 2^{j-1}|I_a(p)\bigr\},
\end{eqnarray*}
which gives (\ref{eq:sieve_4_0}). \qed

\medskip
\noindent It seems that our result on $Q_a(x;4,0)$ cannot be derived from Odoni's result, because $4$ is not square free. So we describe our proof briefly here. 

\medskip
\noindent{\bf Proof of Theorem \ref{th:uncond}.} 

The first term of the right hand side of (\ref{eq:sieve_4_0}) is calculated by the Siegel-Walfisz theorem. As to the other terms, we take
$$\eta_1=\log\log x \quad\mbox{ and }\quad \eta_2=\sqrt x\log x.$$
Then
\begin{eqnarray*}
\displaystyle\sum_{j\geq 2}\sharp
\{p\leq x &;& p\equiv 1\pmod{2^{j}},\ 2^{j-1}|I_a(p)\bigr\}\\
&=&\Bigl(\sum_{2^j\leq\eta_1}+\sum_{\eta_1< 2^j\leq\eta_2}+
\sum_{\eta_2< 2^j\leq x}\Bigr)\sharp
\{p\leq x\ ;\ p\equiv 1\pmod{2^{j}},\ 2^{j-1}|I_a(p)\bigr\}\\
&=&I_1+I_2+I_3, \mbox{ say.}
\end{eqnarray*}
By the Siegel-Walfisz Theorem again, we have
$$I_2\ll\frac{x}{\log x\log\log x},$$
and in a similar way to Hooley \cite{Ho}, we have
$$I_3\ll\frac{x}{\log^3 x}.$$
For a prime $p$, ``$p\equiv 1\pmod{2^i}$ and $2^j|I_a(p)$'' if and only if $p$ splits completely in the field $K_{i,j}^{(2)}$. Thus
$$I_1=\sum_{2^j\leq\eta_1}\frac{1}{[K_{j,j-1}^{(2)}:{\bf Q}]}\pi^{(1)}(x;K_{j,j-1}^{(2)}),$$
and Theorem \ref{th:PIT} gives 
$$I_1=\sum_{j=2}^{\infty}\frac{1}{n_j}\li x-\sum_{2^j>\eta_1}\frac{1}{n_j}\li x
+O\left(\sum_{2^j\leq \eta_1}xe^{-c\sqrt{\log x}/n_j^2}\right),$$
where $n_j=[K_{j,j-1}^{(2)}:{\bf Q}]$. When $2^j\leq \eta_1$, 
$$n_j^2\ll 2^{4j}\ll(\log\log x)^4,$$
and we have
$$\sum_{2^j>\eta_1}\frac{1}{n_j}\li x\ll \frac{x}{\log x(\log\log x)^2},$$
$$\sum_{2^j\leq \eta_1}xe^{-c\sqrt{\log x}/n_j^2}\ll \frac{x}{\log^2 x}.$$
Consequently, (\ref{eq:sieve_4_0}) turn into
\begin{eqnarray}
\sharp Q_a(x;4,0)&=&\left\{\frac{1}{\varphi(4)}-
\sum_{j\geq 1}\left(\frac{1}{[K_{j+1,j}^{(2)}:{\bf Q}]}-\frac{1}{[K_{j+2,j}^{(2)}:{\bf Q}]}\right)\right\}\li x\nonumber\\
&+&O\left(\frac{x}{\log x\log\log x}\right).
\label{eq:degree_4_2}
\end{eqnarray}
Now, using Lemma \ref{th:moree_lemma}, we can easily verify that the coefficient of $\li x$ is equal to $1/3$. 

When $l=2$, we notice that 
$$\sharp Q_a(x;4,2)=\sharp Q_a(x;2,0)-\sharp Q_a(x;4,0).$$
We already have the asymptotic formula for $\sharp Q_a(x;4,0)$, and from Odoni's result, we have
$$\sharp Q_a(x;2,0)=\frac{2}{3}\li x +O\left(\frac{x}{\log x\log\log x}\right).$$
This completes the proof of Theorem \ref{th:uncond}. (Of course we can derive the same result from (\ref{eq:sieve_4_2}) directly.) \qed

\medskip
\rem It is clear from the proof described above that the assumptions $a$ is square free and $a\ne 2$ are not essential. In fact, if we spare no effort to calculate the degrees in (\ref{eq:degree_4_2}) (and in the corresponding formula for $\sharp Q_a(x;2,0)$), we can find the densities $\sharp Q_a(x;4,0)$ and $\sharp Q_a(x;4,2)$ for other types of $a$'s. 

\section{Proof of Theorem \ref{th:cond}}
We shall prove here Theorem \ref{th:cond}. Our proof consists of two parts:

\medskip
\noindent{\bf Part I. } In Proposition \ref{th:4-5} , under the assumptions of Theorem \ref{th:cond}, we prove the existence of the natural densities
$$\delta_j=\lim_{x\to\infty}\frac{\sharp Q_a(x;4,j)}{\pi(x)}$$
for $j=1,3$. 

\medskip
\noindent{\bf Part II. } We prove in Proposition \ref{th:4-7} that 
\begin{equation}\label{eq:4-1}
\delta_1=\delta_3.
\end{equation}
Then our unconditional result Theorem \ref{th:uncond} shows that 
$$\delta_1+\delta_3=1-(\delta_0+\delta_2)=\frac13,$$
and (\ref{eq:4-1}) proves Theorem \ref{th:cond}. 

\medskip
\noindent{\bf Part I. } 

\medskip
We start our proof from formulas (\ref{eq:sieve_4_1}) and (\ref{eq:sieve_4_3}).

\bigskip
\noindent{\bf 1$^\circ$ A reduction of the double-infinite-sum in (\ref{eq:sieve_4_1}).}

\medskip
In order to simplify the double-infinite-sum in (\ref{eq:sieve_4_1}), we apply Lemma \ref{th:lemma_I_a(p)}. We take $\psi(x)=\log\log x$, then
\begin{eqnarray*}
\sum_{2^f+l\cdot 2^{f+2}\geq\log\log x}
\sharp N_a(x;2^f+l\cdot 2^{f+2};1+2^f\pmod{2^{f+2}})
&\leq&\sharp \bigl\{p\leq x\ ;\ I_a(p)\geq\log\log x\bigr\}\\
&\ll&\pi(x)(\log\log x)^{-1}.
\end{eqnarray*}
So we have 
\begin{eqnarray}
\sharp Q_a(x;4,1)
&=&\sum_{f\geq1,l\geq0\atop 2^f+l\cdot 2^{f+2}\leq\log\log x}\sharp 
  N_a(x;2^f+l\cdot 2^{f+2};1+2^f\pmod{2^{f+2}})\nonumber\\
&+&\sum_{f\geq1,l\geq0\atop 2^f+l\cdot 2^{f+2}\leq\log\log x}\sharp 
  N_a(x;3\cdot 2^f+l\cdot 2^{f+2};1+3\cdot 2^f\pmod{2^{f+2}})\nonumber\\
&+&O\left(\pi(x)\frac{1}{\log\log x}\right).
\label{eq:4-3}
\end{eqnarray}
This formula shows that, in our proof, the calculation of $\sharp N_a(x;2^f+l\cdot 2^{f+2};1+2^f\pmod{2^{f+2}})$ is very important. In the literature, $\sharp N_a(x;n)$ is already calculated by Lenstra \cite{Le} and Murata \cite{Mu} for arbitrary $n\in{\bf N}$. And $\{p\leq x\ ;\ I_a(p)=1,\ p\equiv s\pmod t\}$ is already considered in Lenstra \cite{Le} (see also Moree \cite{Mo}). So, in what follows, we combine Lenstra's idea about the control of residue classes and Murata's method of obtaining an asymptotic formula of $\sharp N_a(x;n)$. 

\bigskip
\noindent{\bf 2$^\circ$ A decomposition of $\sharp N_a(x;2^f+l\cdot 2^{f+2};1+2^f\pmod{2^{f+2}})$.}

\medskip
In the formulas (\ref{eq:sieve_4_1}) and (\ref{eq:sieve_4_3}), we find four terms of the same type: $\sharp N_a(x;j\cdot 2^f+l\cdot 2^{f+2};1+j'\cdot 2^f\pmod{2^{f+2}})$, for $(j,j')=(1,1),(1,3),(3,1)$ and $(3,3)$. We can calculate these terms in the same way, and in 2$^\circ$ -- 5$^\circ$ we think about only the case $(j,j')=(1,1)$. 

Here we need some new notations. The letter $q$ always means a prime number, and we put 
\begin{eqnarray*}
k&=&2^f+l\cdot 2^{f+2},\\
k_0&=&\prod_{q|k}q \quad\mbox{(i.e. the core of }k).
\end{eqnarray*}
We define the algebraic number field 
$$K_k={\bf Q}(\zeta_{k_0},a^{1/k})$$
and define two sets of prime ideals of $K_k$:

$$B(x;K_k;a^{1/k};N)
=\left\{\begin{array}{l}
  \mathfrak p: \mbox{ a prime ideal in }K_k,\ N\mathfrak p=p^1\leq x,
     p\equiv 1\pmod N\\
  a^{1/k}\mbox{ is a primitive root }\mod\mathfrak p \end{array}
\right\}$$
and 
$$B(x;K_k;a^{1/k};N;s\pmod t)
=\{\mathfrak p\in B(x;K_k;a^{1/k};N)\ ;\ p\equiv s\pmod t\}.$$

First we decompose $\sharp N_a(x;k;1+2^f\pmod{2^{f+2}})$ into the sum of some $\sharp B(x;K_k;a^{1/k};N;$ $s \pmod t)$'s and here we use Murata's method \cite{Mu}. 

As a consequence of \cite[Lemma 3]{Mu}, we have, if $\mathfrak p\in B(x;K_k;a^{1/k};N;s\pmod t)$, then $N\mathfrak p=p\in N_a(x;k;s\pmod t)$. And the same argument with \cite[Lemma 4]{Mu}, we get easily that, if $p\in N_a(x;k;s\pmod t)$, then $p$ gives rise to $[K_k:{\bf Q}]\varphi(\tilde p)/\tilde p$ elements of $B(x;K_k;a^{1/k};k;s\pmod t)$, where $\tilde p$ is defined in \cite[p. 559]{Mu}. In fact, we can prove these results only by limiting the proofs of \cite{Mu} into the residue class $s\pmod t$. 

From these relations, making use of the M\"obius inversion formula, we can deduce the following decomposition:
\begin{prop}\footnote{We can simplify \cite[Proposition 1]{Mu} into
$$|N_a^{(n)}(x)|=\frac{1}{[K_n:{\bf Q}]}\frac{n_0}{\varphi(n_0)}\sum_{d|n_0}\frac{\mu(d)}{d}|B(\sqrt[n]{a};K_n;x;nd)|.$$ This simplification is due to Dr. R. Takeuchi.}\label{th:4-1}
 \ \\

$\sharp N_a(x;k;1+2^f\pmod{2^{f+2}})$
$$=\frac{1}{[K_k:{\bf Q}]}
\frac{k_0}{\varphi(k_0)}\sum_{d|k_0}\frac{\mu(d)}{d}
\sharp B(x;K_k;a^{1/k};kd;1+2^f\pmod{2^{f+2}}).$$
\end{prop}

\bigskip
\noindent{\bf 3$^\circ$ Calculation of $\sharp B(x;K_k;a^{1/k};kd;1+2^f\pmod{2^{f+2}})$.}

\medskip
The cardinality $\sharp B(x;K_k;a^{1/k};m)$ is already calculated in \cite[Proposition 2]{Mu}, and this calculation is carried out along Hooley's work \cite{Ho}. In our present case, we have to take into account the condition ``$p\equiv 1+2^f\pmod{2^{f+2}}$'', but the calculation itself needs only some slight modifications. We omit the detail and present our result in Proposition \ref{th:4-2}. 

We define
$$P(x;K_k;a^{1/k};kd;n)
=\left\{\begin{array}{l}
  \mathfrak p: \mbox{ a prime ideal in }K_k,\ N\mathfrak p=p^1\leq x,\ 
  p\equiv 1\pmod{kd},\\
  \mbox{the equation }X^q\equiv a^{1/k}\pmod{\mathfrak p} 
  \mbox{ is solvable in }O_{K_k} \mbox{ for any }q|n \end{array}
\right\},$$
where $O_{K_k}$ is the ring of integers of the field $K_k$, and 
$$P(x;K_k;a^{1/k};kd;s\pmod t;n)=
\{\mathfrak p\in P(x;K_k;a^{1/k};kd;n)\ ;\ p\equiv s\pmod t\}.$$
Then we have
\begin{prop}\label{th:4-2}
\ \\

$\sharp B(x;K_k;a^{1/k};kd;1+2^f\pmod{2^{f+2}})$
\begin{equation}\label{eq:4-5}
=\psum_{n}\mu(n)\sharp P(x;K_k;a^{1/k};kd;1+2^f\pmod{2^{f+2}};n)
+O\left(\frac{x(\log\log x)^3}{\log^2 x}\right),
\end{equation}
where the $\psum_n$ means the sum over such an $n\leq x$ which is either $1$ or a positive square free integer composed entirely of prime factors not exceeding $(1/8)\log x$, and the constant implied by the $O$-symbol is absolute (In Hooley's paper \cite{Ho}, he made use of $(1/6)\log x$ instead of $(1/8)\log x$). 
\end{prop}

\bigskip
\noindent{\bf 4$^\circ$ Calculation of $\sharp P(x;K_k;a^{1/k};kd;1+2^f\pmod{2^{f+2}};n)$.}

\medskip
Here we need GRH. 

We define algebraic extension fields
\begin{eqnarray*}
G_{k,n,d}&=&K_k(\zeta_n,a^{1/kn},\zeta_{kd}),\\
\tilde G_{k,n,d}&=&G_{k,n,d}(\zeta_{2^{f+2}}),
\end{eqnarray*}
and we take an automorphism $\sigma\in\Aut({\bf Q}(\zeta_{2^{f+2}})/{\bf Q})$ which is defined by 
$$\sigma:\zeta_{2^{f+2}} \mapsto (\zeta_{2^{f+2}})^{1+2^f}.
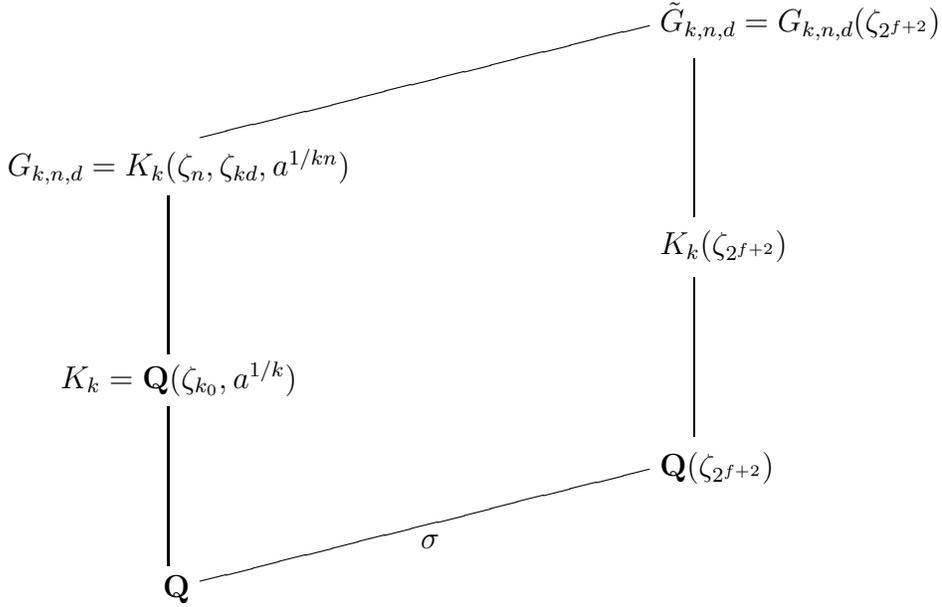
$$
\begin{center}
\begin{picture}(270,250)
\put(4,0){\makebox(20,14)[l]{${\bf Q}$}}
\put(0,80){\makebox(20,14)[c]{$K_k={\bf Q}(\zeta_{k_0}, a^{1/k})$}}
\put(0,160){\makebox(20,14)[c]{$G_{k,n,d}=K_k(\zeta_n, \zeta_{kd}, a^{1/kn})$}}
\put(192,46){\makebox(20,14)[l]{${\bf Q}(\zeta_{2^{f+2}})$}}
\put(192,130){\makebox(20,14)[l]{$K_k(\zeta_{2^{f+2}})$}}
\put(192,215){\makebox(20,14)[l]{$\tilde G_{k,n,d}=G_{k,n,d}(\zeta_{2^{f+2}})$}}
%
\put(6,16){\line(0,1){60}}
\put(6,96){\line(0,1){60}}
\put(205,65){\line(0,1){60}}
\put(205,148){\line(0,1){60}}
%
\put(18,10){\line(4,1){170}}
\put(18,178){\line(4,1){170}}
%
\put(95,18){\makebox(20,14){$\sigma$}}
\end{picture}
\end{center}
Let $\sigma^\ast\in\Aut(\tilde G_{k,n,d}/K_k)$ be the automorphism defined by 
\begin{equation}\label{eq:4-6}
\left\{\begin{array}{l}
\sigma^\ast|_{G_{k,n,d}}=\id_{G_{k,n,d}},\\[2mm]
\sigma^\ast|_{{\bf Q}(\zeta_{2^{f+2}})}=\sigma.
\end{array}\right.
\end{equation}
Here we remark that such a $\sigma^\ast$ does not always exist. When we can construct this $\sigma^\ast$ from $\sigma$, we can prove the following Lemma \ref{th:4-3} and Proposition \ref{th:4-4}:
\begin{lem}\label{th:4-3}
When $\sigma^\ast$ exists, $\{\sigma^\ast\}$ is a conjugacy class of $\Aut(\tilde G_{k,n,d}/K_k)$ by itself. 
\end{lem}
\begin{prop}\label{th:4-4}
We assume GRH, and $\sigma^\ast\in\Aut(\tilde G_{k,n,d}/K_k)$ exists. Then
\begin{equation}
\sharp P(x;K_k;a^{1/k};kd;1+2^f\pmod{2^{f+2}};n)=
\pi(x;\tilde G_{k,n,d}/K_k,\{\sigma^\ast\})+O(k^2\sqrt{x}(\log\log x)^4),
\label{eq:4-7}
\end{equation}
where $\pi(x;\tilde G_{k,n,d}/K_k,\{\sigma^\ast\})$ is defined in (\ref{eq:def_pi_C}). 
\end{prop}
When the automorphism $\sigma^\ast\in\Aut(\tilde G_{k,n,d}/K_k)$ does not exist, 
we regard $\pi(x;\tilde G_{k,n,d}/K_k,\{\sigma^\ast\})=0$. 

\medskip
\noindent{\bf Proof of Lemma \ref{th:4-3}. } We take an arbitrary $\tau\in\Aut(\tilde G_{k,n,d}/K_k)$. Since $\tau(a^{1/nk})=a^{1/nk}{\zeta_n}^i$ for some $i\in{\bf N}$, we have 
$$\sigma^\ast\circ\tau(a^{1/nk})=a^{1/nk}{\zeta_n}^i
=\tau\circ\sigma^\ast(a^{1/nk}).$$
Similarly we can prove 
$$\sigma^\ast\circ\tau(\zeta_n)=\tau\circ\sigma^\ast(\zeta_n)\quad \mbox{ and }
\quad \sigma^\ast\circ\tau(\zeta_{kd})=\tau\circ\sigma^\ast(\zeta_{kd}).$$
Moreover, since $\tau(\zeta_{2^{f+2}})=(\zeta_{2^{f+2}})^{i'}$ for some $i'\in{\bf N}$, we have
$$\sigma^\ast\circ\tau(\zeta_{2^{f+2}})=(\zeta_{2^{f+2}})^{i'(1+2^f)}
=\tau\circ\sigma^\ast(\zeta_{2^{f+2}}).$$
These prove our assertion. \qed

\medskip
\noindent{\bf Proof of Proposition \ref{th:4-4}. } We notice that the next two conditions are equivalent:

\medskip
(a) $\mathfrak p\in P(x;K_k;a^{1/k};kd;n)$,

(b) $\mathfrak p$ splits completely in the extension $G_{k,n,d}/K_k$ and $N\mathfrak p\leq x$, 

\medskip
\noindent and also the following two are equivalent:

\medskip
(c) $p\equiv 1+2^f\pmod{2^{f+2}}$,

(d) the Frobenius map $(p,{\bf Q}(\zeta_{2^{f+2}})/{\bf Q})=\sigma$.

\medskip
\noindent Now let $\mathfrak p\in P(x;K_k;a^{1/k};kd;1+2^f\pmod{2^{f+2}};n)$. $\tilde G_{k,n,d}/K_k$ is a normal extension, and we can define $[\mathfrak p,\tilde G_{k,n,d}/K_k]$, the conjugacy class of Frobenius automorphisms corresponding to prime ideals $\mathfrak P\subset \tilde G_{k,n,d}$ over $\mathfrak p$. We take $\tau\in[\mathfrak p,\tilde G_{k,n,d}/K_k]$. Then the equivalent relation (a) $\Leftrightarrow$ (b) implies that the ideal $\mathfrak p$ splits completely in the extension $G_{k,n,d}/K_k$, thus $\tau|_{G_{k,n,d}}=\id_{G_{k,n,d}}$. Also the equivalent relation (c) $\Leftrightarrow$ (d) implies $\tau|_{{\bf Q}(\zeta_{2^{f+2}})}=\sigma$. Consequently $\tau=\sigma^\ast$ and , with Lemma \ref{th:4-3}, $[\mathfrak p,\tilde G_{k,n,d}/K_k]=\{\sigma^\ast\}$. This proves

\medskip
$P(x;K_k;a^{1/k};kd;1+2^f\pmod{2^{f+2}};n)$
$$\subset
\{\mathfrak p:\mbox{ a prime ideal in }K_k,[\mathfrak p,\tilde G_{k,n,d}/K_k]=\{\sigma^\ast\}, N\mathfrak p\leq x\}.$$
Now we consider the prime ideal $\mathfrak p$ of $K_k$ which satisfies two conditions: 
\begin{equation}
[\mathfrak p,\tilde G_{k,n,d}/K_k]=\{\sigma^\ast\}\quad\mbox{ and }\quad
N\mathfrak p\leq x.
\label{eq:4-8}
\end{equation}
It is easy to see that
$$\sharp \left\{\begin{array}{l}
\mathfrak p: \mbox{ a prime ideal of }K_k, 
\mathfrak p \mbox{ satisfies the conditions (\ref{eq:4-8})},\\
N\mathfrak p=p^s \mbox{ with }s\geq2
\end{array}\right\}$$
\begin{eqnarray*}
&\ll&[K_k:{\bf Q}]\cdot O\left(\sum_{i=2}^\infty x^{1/i}\right)\\
&\ll& k^2\sqrt x\bigl(\log\log x\bigr)^4.
\end{eqnarray*}
This means that, except for at most $O\bigl(k^2\sqrt x\bigl(\log\log x\bigr)^4\bigr)$ of primes, we can assume $N\mathfrak p=p\leq x$. Then the property $\sigma^\ast|_{G_{k,n,d}}=\id_{G_{k,n,d}}$ implies that $\mathfrak p$ satisfies the condition (b), and $\mathfrak p\in P(x;K_k;a^{1/k};kd;n)$. Furthermore, the property $\sigma^\ast|_{{\bf Q}(\zeta_{2^{f+2}})}=\sigma$ implies (d). These show that $\mathfrak p\in P(x;K_k;a^{1/k};kd;1+2^f\pmod{2^{f+2}};n)$, and we proved our assertion. \qed

\bigskip
\noindent{\bf 5$^\circ$ Existence of the densities $\delta_1$ and $\delta_3$.}

\medskip
We can now prove the main result of Part I. 
\begin{prop}\label{th:4-5}
We assume GRH. Then for $j=1,3$, $\sharp Q_a(x;4,j)$ has the natural density $\delta_j$. 
\end{prop}
Before we prove this proposition, we prepare some estimates:
\begin{lem}\label{th:4-6}
Under the above notations, let $d_{\tilde G_{k,n,d}}$ be the discriminant of the field $\tilde G_{k,n,d}$. Then

\medskip
\noindent{\rm (i)}
$$[\tilde G_{k,n,d}:K_k]=\delta\frac{d}{k_0\varphi((n,k_0))}\cdot kn\varphi(n),$$
where $\delta$ is one of the five numbers $\{8,4,2,1,1/2\}$. 

\medskip
\noindent{\rm (ii)}
$$\log|d_{\tilde G_{k,n,d}}|\ll(nkd)^3\log(nkd).$$
\end{lem}
\prf It is already proved in Murata \cite{Mu} that $[K_k:{\bf Q}]=\eta_1k\varphi(k_0)$, where $\eta_1=1$ or $1/2$, and that ({\it cf.} formulas (11) and (12) of \cite{Mu})
$$[G_{k,n,d}:{\bf Q}]=\eta_2\frac{\varphi(k_0)k^2d}{k_0}\frac{n\varphi(n)}{\varphi((n,k_0))},$$
where $\eta_2=1$ or $1/2$. Moreover, since $\tilde G_{k,n,d}=G_{k,n,d}(\zeta_{2^{f+2}})$ and $\zeta_{2^f}\in G_{k,n,d}$, 
$$[\tilde G_{k,n,d}:G_{k,n,d}]\Bigm|4.$$
Combining these formulas, we get (i) easily. 

We now prove (ii). Let $L_1/{\bf Q}$ and $L_2/{\bf Q}$ be two extension fields, $L$ be the composite field $L_1\cdot L_2$, and $d_{L_1}$, $d_{L_2}$, $d_L$ be the discriminants of $L_1$, $L_2$, $L$, respectively. Then we have the following relation:
$$|d_L|\Bigm| |d_{L_1}|^{[L:L_1]}|d_{L_2}|^{[L:L_2]}.$$
From this, we have an estimate
\begin{equation}\label{eq:4-9}
|d_L|\leq |d_{L_1}|^{[L_2:{\bf Q}]}|d_{L_2}|^{[L_1:{\bf Q}]}.
\end{equation}
Here we take 
\begin{eqnarray*}
L_1&=&{\bf Q}(a^{1/nk}),\\
L_2&=&{\bf Q}(\zeta_n,\zeta_{kd}, \zeta_{2^{f+2}}).
\end{eqnarray*}
It is known that the discriminant of the cyclotomic field ${\bf Q}(\zeta_{p^r})$ is given by 
$$|d_{{\bf Q}(\zeta_{p^r})}|=p^{p^{r-1}(pr-r-1)}.$$
From this, it is easy to prove that, for any $m\in{\bf N}$, 
$$\log|d_{{\bf Q}(\zeta_{m})}|\leq m^2\log m.$$
Thus we have firstly,
\begin{equation}\label{eq:4-10}
\log|{d_{L_2}}|\leq(4nkd)^2\log(4nkd).
\end{equation}
We also have
\begin{eqnarray*}
|d_{L_1}|&\leq&\Bigl|\mbox{the discriminant of the polynomial }X^{nk}-a\Bigr|\\
&=&\left|\prod_{0<i<j<nk}a^{2/nk}({\zeta_{nk}}^i-{\zeta_{nk}}^j)^2\right|\\
&\leq& a^{nk}|d_{{\bf Q}(\zeta_{nk})}|,
\end{eqnarray*}
and then
\begin{equation}\label{eq:4-11}
\log|{d_{L_1}}|\leq nk\log a+(nk)^2\log (nk).
\end{equation}
Since $[L_1:{\bf Q}]\leq nk$ and $[L_2:{\bf Q}]\leq 4nkd$, we now prove (ii) from (\ref{eq:4-9}), (\ref{eq:4-10}) and (\ref{eq:4-11}). \qed
\begin{cor}\label{th:4-7}
We assume GRH. The numbers $k,n,d$ are as above, and let $k\leq\log\log x$. Then, for any $\varepsilon>0$, we have
\begin{equation}\label{eq:4-12}
\pi(x;\tilde G_{k,n,d}/K_k, \{\sigma^\ast\})
=\frac{1}{[\tilde G_{k,n,d}:K_k]}\li x
+O\left(x^{\frac34+2\varepsilon}(\log x)^2\right),
\end{equation}
where the constant involved by the $O$-symbol depends only on $\varepsilon$. 
\end{cor}

\medskip
\prf Lemma \ref{th:4-3} says that, when $\sigma^\ast$ exists, we can take the conjugacy class $C=\{\sigma^\ast\}$ and $\sharp C=1$. Let us apply Theorem \ref{th:Chebotarev} for 
$$L=\tilde G_{k,n,d},\ K=K_k,\ C=\{\sigma^\ast\},$$
then, in order to prove Corollary \ref{th:4-7}, it is now sufficient to estimate the next two terms:
$$\frac{\sqrt x}{[\tilde G_{k,n,d}:K_k]}\log d_{\tilde G_{k,n,d}}\quad
\mbox{ and }\quad \log d_{\tilde G_{k,n,d}}.$$
Here we recall that $n$ is a square free integer composed entirely of prime factors not exceeding $(1/8)\log x$. Then we can estimate
$$\log n\leq \sum_{p\leq \frac18\log x}\log p\ll \left(\frac{1}{8}+\varepsilon\right)\log x,$$
and we have 
$$n\ll x^{\frac{1}{8}+\varepsilon}$$
for any $\varepsilon>0$. Then by Lemma \ref{th:4-6} (i) and (ii), 
\begin{eqnarray*}
\frac{\sqrt x}{[\tilde G_{k,n,d}:K_k]}\log d_{\tilde G_{k,n,d}}&\ll&
 \sqrt xk_0(nkd)^2\log(nkd)\\
&\ll& x^{\frac34+2\varepsilon}\log x(\log\log x)^5,
\end{eqnarray*}
$$\log d_{\tilde G_{k,n,d}}\ll x^{\frac38+3\varepsilon}\log x(\log\log x)^6,$$
and these prove (\ref{eq:4-12}). \qed

\medskip
\noindent{\bf Proof of Proposition \ref{th:4-5}. } We define here the number $c(n)$ by 
$$c(n)=\left\{\begin{array}{ll}
1, & \mbox{if the Frobenius map }\sigma^\ast \mbox{ defined by (\ref{eq:4-6})
   exists},\\
0, & \mbox{if not.}
\end{array}\right.$$
Combining (\ref{eq:4-5}), (\ref{eq:4-7}) and Corollary \ref{th:4-7}, we have, provided $k\leq\log\log x$, 

\bigskip
$\sharp B(x;K_k;a^{1/k};kd;1+2^f\pmod{2^{f+2}})$
\begin{eqnarray*}
&=&\psum_{n}\mu(n)c(n)\Bigl(\pi(x;\tilde G_{k,n,d}/K_k,\{\sigma^\ast\})
  +O(\sqrt x(\log\log x)^6)\Bigr)
  +O\left(\frac{x(\log\log x)^3}{\log^2x}\right)\\
&=&\psum_{n}\frac{\mu(n)c(n)}{[\tilde G_{k,n,d}:K_k]}\li x
  +O\left(x^{\frac78+3\varepsilon}(\log x)^2\right)
  +O\left(\frac{x(\log\log x)^3}{\log^2x}\right).
\end{eqnarray*}
Now, from Lemma \ref{th:4-6} (i), it is seen that the leading coefficient is an absolutely convergent series. Making use of an estimate
$$\psum_{n}\frac{\mu(n)c(n)}{[\tilde G_{k,n,d}:K_k]}
=\sum_{n=1}^\infty\frac{\mu(n)c(n)}{[\tilde G_{k,n,d}:K_k]}
+O\left(\frac{1}{\log x}\right),$$
we obtain a formula
$$\sharp B(x;K_k;a^{1/k};kd;1+2^f\pmod{2^{f+2}})=\tilde\delta_{k,d}\li x
  +O\left(\frac{x(\log\log x)^3}{\log^2x}\right),$$
where 
$$\tilde\delta_{k,d}=\sum_{n=1}^\infty\frac{\mu(n)c(n)}{[\tilde G_{k,n,d}:K_k]}.$$
Then Proposition \ref{th:4-1} gives
\begin{eqnarray*}
\sharp N_a(x;k;1+2^f\pmod{2^{f+2}})&=&\left(\frac{1}{[K_k:{\bf Q}]}
  \frac{k_0}{\varphi(k_0)}
  \sum_{d|k_0}\frac{\mu(d)}{d}\tilde\delta_{k,d}\right)
  \li x\nonumber\\
& &+O\left(\frac{k_0}{k\varphi(k_0)^2}\sum_{d|k_0}\frac{1}{d}
  \frac{x(\log\log x)^3}{\log^2x}\right)\nonumber\\
&=&\tilde\delta_{k}\li x+O\left(\frac{x(\log\log x)^3}{\log^2x}\right),
\end{eqnarray*}
with 
$$\tilde\delta_{k}=\frac{1}{[K_k:{\bf Q}]}\frac{k_0}{\varphi(k_0)}
\sum_{d|k_0}\frac{\mu(d)}{d}\tilde\delta_{k,d}.$$
Similarly for $m=3\cdot 2^l+l\cdot 2^{f+2}$, we have
\begin{eqnarray*}
\sharp N_a(x;m;1+3\cdot 2^f\pmod{2^{f+2}})
&=&\tilde\delta_{m}\li x+O\left(\frac{x(\log\log x)^3}{\log^2x}\right).
\end{eqnarray*}
Then, our formula (\ref{eq:4-3}) yields
\begin{eqnarray}
\sharp Q_a(x;4,1)
&=&\left(\sum_{k\leq\log\log x}\tilde\delta_{k}\right)\li x\nonumber\\
&+&\left(\sum_{m\leq\log\log x}\tilde\delta_{m}\right)\li x
  +O\left(\frac{x}{\log x \log\log x}\right).
\label{eq:4-15}
\end{eqnarray}
From the definition, $\tilde\delta_{k}$ and $\tilde\delta_{m}$ are non-negative numbers, and a priori, $\sum_{k=1}^{\infty}\tilde\delta_{k}\leq 1$ and $\sum_{m=1}^{\infty}\tilde\delta_{m}\leq 1$. Thus the two leading coefficients which appeared in (\ref{eq:4-15}) converge, namely,
$$\sum_{k\leq\log\log x}\tilde\delta_{k}+
\sum_{m\leq\log\log x}\tilde\delta_{m}
=\delta_1+o(1).$$
This proves the existence of the density $\delta_1$, and similarly, we can show the existence of $\delta_3$. \qed

\bigskip
\noindent{\bf Part II.}

In this paragraph, we present our proof for $\delta_1=\delta_3$. 

Our proof is based on the following expressions (\ref{eq:4-16}) and (\ref{eq:4-16+}) which we obtained in Part I, but prior to the details, we need a few new notations:
\begin{eqnarray*}
k&=&k(l,f)=2^f+l\cdot 2^{f+2},\\
m&=&m(l,f)=3\cdot 2^f+l\cdot 2^{f+2},
\end{eqnarray*}
$$k_0=\prod_{q|k}q,\quad m_0=\prod_{q|m}q \quad
(\mbox{the cores of $k$ and $m$}),$$
and accordingly
$$\begin{array}{ll}
\tilde G_{k,n,d}=K_k(\zeta_n,a^{1/nk},\zeta_{nd},\zeta_{2^{f+2}}), & 
(\mbox{the same as in Part I})\\
\tilde G_{m,n,d}=K_m(\zeta_n,a^{1/nm},\zeta_{nd},\zeta_{2^{f+2}}). & 
\end{array}$$
Furthermore,
\begin{eqnarray*}
c_1(k,n,d)&=&\left\{\begin{array}{ll}
  1, & \mbox{if we can construct }\sigma_1^\ast\in\Aut(\tilde G_{k,n,d}/K_k)\\
     & \mbox{with the properties (\ref{eq:4-6}),}\\
  0, & \mbox{if not,}
  \end{array}\right.\\
c_3(k,n,d)&=&\left\{\begin{array}{ll}
  1, & \mbox{if we can construct }\sigma_3^\ast\in\Aut(\tilde G_{k,n,d}/K_k)\\
     & \mbox{with the properties (\ref{eq:4-6'}),}\\
  0, & \mbox{if not,}
  \end{array}\right.
\end{eqnarray*}
where
\begin{equation}\label{eq:4-6'}
\left\{\begin{array}{l}
\sigma_3^\ast|_{G_{k,n,d}}=\id_{G_{k,n,d}},\\[2mm]
\sigma_3^\ast|_{{\bf Q}(\zeta_{2^{f+2}})}=\sigma': 
\zeta_{2^{f+2}}\mapsto(\zeta_{2^{f+2}})^{1+3\cdot 2^f}.
\end{array}\right.
\end{equation}
Then
\begin{eqnarray}
\delta_1&=&\sum_{f\geq1}\sum_{l\geq0}\frac{1}{[K_k:{\bf Q}]}
   \frac{k_0}{\varphi(k_0)}\sum_{d|k_0}\frac{\mu(d)}{d}\sum_{n}
   \frac{\mu(n)c_1(k,n,d)}{[\tilde G_{k,n,d}:K_k]}\nonumber\\
&+&\sum_{f\geq1}\sum_{l\geq0}\frac{1}{[K_m:{\bf Q}]}
   \frac{m_0}{\varphi(m_0)}\sum_{d|m_0}\frac{\mu(d)}{d}\sum_{n}
   \frac{\mu(n)c_3(m,n,d)}{[\tilde G_{m,n,d}:K_m]},
\label{eq:4-16}\\
\delta_3&=&\sum_{f\geq1}\sum_{l\geq0}\frac{1}{[K_k:{\bf Q}]}
   \frac{k_0}{\varphi(k_0)}\sum_{d|k_0}\frac{\mu(d)}{d}\sum_{n}
   \frac{\mu(n)c_3(k,n,d)}{[\tilde G_{k,n,d}:K_k]}\nonumber\\
&+&\sum_{f\geq1}\sum_{l\geq0}\frac{1}{[K_m:{\bf Q}]}
   \frac{m_0}{\varphi(m_0)}\sum_{d|m_0}\frac{\mu(d)}{d}\sum_{n}
   \frac{\mu(n)c_1(m,n,d)}{[\tilde G_{m,n,d}:K_m]}.
\label{eq:4-16+}
\end{eqnarray}
We remark here that, when 
\begin{equation}\label{eq:4-17}
\begin{array}{c}
c_1(k,n,d)=c_3(k,n,d),\\
c_1(m,n,d)=c_3(m,n,d),
\end{array}
\end{equation}
then the first term of the right hand side of (\ref{eq:4-16}) is equal to the first term of (\ref{eq:4-16+}), and the second terms of (\ref{eq:4-16}) and (\ref{eq:4-16+}) coincide with each other. This means $\delta_1=\delta_3$. 

Now we prove (\ref{eq:4-17}). 
\begin{prop}\label{th:4-8}
We assume $a_1\equiv 1\pmod 4$. Then the relations (\ref{eq:4-17}) hold. 
\end{prop}
\prf Here we give the proof only for the first relation. 

\medskip
\noindent{\bf Case 1. } If $f\geq 2$, then it is easy to see that, when $c_1(k,n,d)=1$, i.e. $\sigma_1^\ast$ exists, then $(\sigma_1^\ast)^3$ satisfies (\ref{eq:4-6'}), i.e. $c_3(k,n,d)=1$, and vice versa. Thus $c_1(k,n,d)=c_3(k,n,d)$. 

\medskip
\noindent{\bf Case 2. } If $d$ is even, then we can prove 
$$c_1(k,n,d)=c_3(k,n,d)=0.$$
In fact, in this case, $2^{f+1}|kd$, and two conditions
\begin{eqnarray*}
\sigma_1^\ast(\zeta_{kd})&=&\zeta_{kd},\\
\sigma_1^\ast(\zeta_{2^{f+2}})&=&(\zeta_{2^{f+2}})^{1+2^f}
\end{eqnarray*}
contradict each other. Thus $c_1(k,n,d)=0$, and similarly $c_3(k,n,d)=0$. 

\medskip
\noindent{\bf Case 3. } Here we assume $f=1$, $d$ is odd and furthermore $a_1\equiv 1\pmod 4$. We calculate the extension degree $[G_{k,n,d}\cap{\bf Q}(\zeta_8):{\bf Q}]$. 
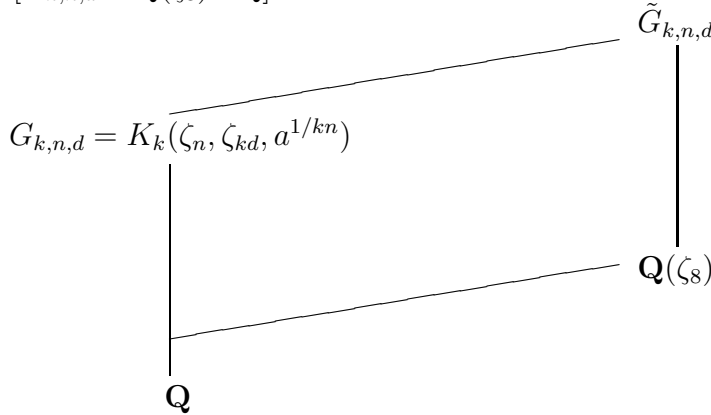
\begin{center}
\begin{picture}(200,150)
\put(4,0){\makebox(20,14)[l]{${\bf Q}$}}
\put(0,100){\makebox(20,14)[c]{$G_{k,n,d}=K_k(\zeta_n, \zeta_{kd}, a^{1/kn})$}}
\put(183,49){\makebox(20,14)[l]{${\bf Q}(\zeta_{8})$}}
\put(183,143){\makebox(20,14)[l]{$\tilde G_{k,n,d}$}}
%
\put(6,16){\line(0,1){80}}
\put(198,65){\line(0,1){76}}
%
\put(6,30){\line(6,1){170}}
\put(6,115){\line(6,1){170}}
\end{picture}
\end{center}
Let $\langle a,b\rangle$ mean the least common multiple of $a$ and $b$. We have 
\begin{eqnarray*}
[G_{k,n,d}:{\bf Q}]&=&\left\{\begin{array}{ll}
nk\varphi(\langle n,kd\rangle), & \mbox{if }2a_1\nmid \langle n,kd\rangle,\\
\Frac12nk\varphi(\langle n,kd\rangle), & \mbox{if }2a_1| \langle n,kd\rangle,
  \end{array}\right.\\
{[\tilde G_{k,n,d}:{\bf Q}]}&=&\left\{\begin{array}{ll}
4nk\varphi(\langle n,kd\rangle), & \mbox{if }2a_1\nmid 4\langle n,kd\rangle,\\
2nk\varphi(\langle n,kd\rangle), & \mbox{if }2a_1|4\langle n,kd\rangle,
  \end{array}\right.
\end{eqnarray*}
and $[{\bf Q}(\zeta_8):{\bf Q}]=4$. Since the two conditions ``$2a_1\nmid \langle n,kd\rangle$'' and ``$2a_1\nmid 4\langle n,kd\rangle$'' are equivalent, we have
$$[G_{k,n,d}\cap {\bf Q}(\zeta_8):{\bf Q}]=
\frac{[G_{k,n,d}:{\bf Q}][{\bf Q}(\zeta_8):{\bf Q}]}
{[\tilde G_{k,n,d}:{\bf Q}]}=1.$$
This means we can construct $\sigma_1^\ast$ with the properties (\ref{eq:4-6}), i.e. $c_1(k,n,d)=1$. Similarly $c_3(k,n,d)=1$. This completes the proof of (\ref{eq:4-17}). \qed

\medskip
As we described in the above, now we can conclude that 
$$\delta_1=\delta_3=\frac16.$$
\section{Numerical Examples}
In this section we look at some numerical calculations of the densities of $Q_a(x;4,l)$, including those for $a$'s which are not dealt with in the previous sections. We did computer calculations of $\sharp Q_a(x;4,l)/\pi(x)$ up to $x=10^7$, where $\pi(x)$ denotes the number of primes not exceeding $x$. 

The tables below exhibit what the densities of $Q_a(x;4,l)$ are like for various square free $a$'s. The exact densities of $Q_a(x;4,0)$ and $Q_a(x;4,2)$ can be found unconditionally by Theorem \ref{th:uncond}. When $a\equiv 1\pmod 4$ (the cases $a=5$ and $21$ in the tables), the exact densities of $Q_a(x;4,1)$ and $Q_a(x;4,3)$ can be proved to be $1/6$ under GRH (Theorem \ref{th:cond}). On the other hand, when $a\not\equiv 1\pmod 4$, the exact densities of $Q_a(x;4,1)$ and $Q_a(x;4,3)$ are unknown even if we assume GRH. Among such $a$'s, in the case $a=3$, both the densities of $Q_3(x;4,1)$ and $Q_3(x;4,3)$ seem very close to $1/6$. But the calculation for $a=6$ shows that there really exists a case when the densities of $Q_a(x;4,1)$ and $Q_a(x;4,3)$ seem different values from $1/6$. This observation shows that the condition $a\equiv 1\pmod 4$ is not just for technical reasons, but plays an essential role in Theorem \ref{th:cond} for determining of the natural densities of $\sharp Q_a(x;4,l)$, $l=1,3$. 
\medskip
\begin{center}
{\small {\bf Table 5.1.} The densities of $Q_{5}(x;4,l)$}
\nopagebreak

\smallskip
\begin{tabular}{|c||c|c|c|c|}\hline
  $x$   &  $l=0$   &  $l=1$   &  $l=2$   &  $l=3$   \\ \hline
$10^3$   & 0.319277 & 0.156627 & 0.349398 & 0.174699\\
$10^4$   & 0.327628 & 0.167074 & 0.340668 & 0.164629\\
$10^5$   & 0.334619 & 0.167049 & 0.333055 & 0.165276\\
$10^6$   & 0.333227 & 0.167155 & 0.332934 & 0.166684\\
$10^7$   & 0.333320 & 0.166771 & 0.333099 & 0.166810\\\hline
\end{tabular}

\medskip
{\small {\bf Table 5.2.} The densities of $Q_{21}(x;4,l)$}
\nopagebreak

\smallskip
\begin{tabular}{|c||c|c|c|c|}\hline
  $x$   &  $l=0$   &  $l=1$   &  $l=2$   &  $l=3$   \\ \hline
$10^3$   & 0.339394 & 0.133333 & 0.339394 & 0.187879\\
$10^4$   & 0.329527 & 0.160685 & 0.334421 & 0.175367\\
$10^5$   & 0.333507 & 0.166649 & 0.333194 & 0.166649\\
$10^6$   & 0.332582 & 0.165972 & 0.334110 & 0.167335\\
$10^7$   & 0.332836 & 0.166527 & 0.333917 & 0.166720\\ \hline
\end{tabular}

\medskip
{\small {\bf Table 5.3.} The densities of $Q_{3}(x;4,l)$}
\nopagebreak

\smallskip
\begin{tabular}{|c||c|c|c|c|}\hline
  $x$   &  $l=0$   &  $l=1$   &  $l=2$   &  $l=3$   \\ \hline
$10^3$   & 0.331325 & 0.150602 & 0.331325 & 0.186747\\
$10^4$   & 0.331703 & 0.163814 & 0.339038 & 0.165444\\
$10^5$   & 0.334411 & 0.167362 & 0.332325 & 0.165902\\
$10^6$   & 0.332488 & 0.166607 & 0.333762 & 0.167142\\
$10^7$   & 0.333298 & 0.166757 & 0.333397 & 0.166548\\ \hline
\end{tabular}

\medskip
{\small {\bf Table 5.4.} The densities of $Q_{6}(x;4,l)$}
\nopagebreak

\smallskip
\begin{tabular}{|c||c|c|c|c|}\hline
   $x$    &  $l=0$   &  $l=1$   &  $l=2$   &  $l=3$\\ \hline
$10^3$ & 0.331325 & 0.126506 & 0.325301 & 0.216867\\
$10^4$ & 0.334963 & 0.133659 & 0.333333 & 0.198044\\
$10^5$ & 0.333785 & 0.133577 & 0.332847 & 0.199791\\
$10^6$ & 0.333151 & 0.132249 & 0.333507 & 0.201093\\
$10^7$ & 0.333331 & 0.132179 & 0.333019 & 0.201471\\ \hline
\end{tabular}
\end{center}
 
\end{document}